%AMSteX Version 2.2
\input amstex
\documentstyle{amsppt}
\magnification=\magstep 1
%\NoBlackBoxes
\voffset-1cm
\TagsOnLeft
%\subjclassyear{2000}
\loadbold
\topmatter
\title Nielsen coincidence theory in arbitrary  codimensions  \endtitle
%\rightheadtext{}
\author Ulrich Koschorke *\endauthor
\leftheadtext{Ulrich Koschorke}
%\affil Universit\"at Siegen, Emmy-Noether-Campus, D-57068 Siegen \endaffil
\address Universit\"at Siegen,
Emmy Noether Campus, Walter-Flex-Str. 3,
D-57068 Siegen, Germany
\endaddress
\email koschorke\@mathematik.uni-siegen.de \endemail
\thanks * Supported in part within the German-Brazilian Cooperation by IB-BMBF.
\endthanks
\abstract
Given two maps \ $f_1, f_2 : M^m \longrightarrow N^n$ \ between manifolds of the indicated arbitrary dimensions, when can they be deformed away from one another? More generally: what is the minimum number \ $MCC (f_1, f_2)$ \ of {\it pathcomponents} \  of the coincidence space of maps \ $f'_1, f'_2$ \ where \ $f'_i$ \ is homotopic to \ $f_i, \ i = 1, 2\ $? Approaching this question via normal bordism theory we define a lower bound \ $N (f_1, f_2)$ \ which generalizes the Nielsen number studied in classical fixed point and coincidence theory (where \ $m = n$). In at least three settings \ $N (f_1, f_2)$ \ turns out to coincide with \ $MCC (f_1, f_2)$: \ (i) \ when \ $m < 2n - 2$; (ii) when \ $N$ is the unit circle; and (iii) when \ $M$ \ and \ $N$ \ are spheres and a certain injectivity condition involving James-Hopf invariants is satisfied. We also exhibit situations where \ $N (f_1, f_2)$ \ vanishes, but \ $MCC (f_1, f_2)$ \ is strictly positive.
\endabstract
\subjclass
Primary 55 M 20, 55 S 35, 57 R 90. Secondary 55 P 35, 55 Q 25, 55 Q 40, 57 R 42
\endsubjclass
\endtopmatter

%%%%% eingefuegt fr Bilder und Textrotation
% \input graphtex.tex
\input BoxedEPS.tex
\SetRokickiEPSFSpecial
\HideDisplacementBoxes
%%%%% ende %%%%%%%%%%%%%%

\define\incl{\operatorname{incl}}
%\define\fill{\operatorname{fill}}
%\define\quot{\operatorname{quot}}
%\define\sign{\operatorname{sign}}
%\define\coll{\operatorname{coll}}
\define\rel{\operatorname{rel}}
\define\coker{\operatorname{coker}}
\define\Iso{\operatorname{Iso}}
\define\id{\operatorname{id}}
\define\inv{\operatorname{inv}}

\define\scirc{{\ssize{\circ}}}
\document

\specialhead 1.\ \ Introduction and statement of results
\endspecialhead

Throughout this paper \ $f_1, f_2 : M \longrightarrow N$ \ denote two (continuous) maps between given smooth connected manifolds \ $M$  \ and \ $N$ \ without boundary, of strictly positive dimensions \ $m$ \ and \ $n$, \ resp., \ $M$ \ being compact. (All manifolds are assumed Hausdorff with a countable basis).

\definition{Definition 1.1} \ We say the pair \ $(f_1, f_2)$ \ is \ {\it loose}\ \ if there exist (continuous) maps \ $f'_i$ \ homotopic to \ $f_i, \ i = 1, 2$, \ without coincidences (i.e. $f'_1 (x) \ne f'_2 (x)$ \ for all $x \in M)$ \ or, equivalently, if the {\it {\bf m}inimum number of {\bf c}oincidence {\bf c}omponents} \
$$
MCC (f_1, f_2) \ := \ \min \{ \# \pi_0 (C (f'_1, f'_2)) | f'_1 \sim f_1, f'_2 \sim f_2\}
$$
vanishes. Here \ $\# \pi_0 (C (f'_1, f'_2))$ \ denotes the number of {\it pathcomponents} \  of the coincidence subspace
$$
C (f'_1, f'_2) \ = \ \{ x \in M | f'_1 (x) = f'_2 (x) \}
$$
of \ $M$.
\enddefinition

\subhead Question \endsubhead \ When is \ $(f_1, f_2)$ \ loose? \ In other words: when can the maps \ $f_1$ \ and \ $f_2$ \ be deformed away from one another? More generally, how big is \ $MCC (f_1, f_2)$?
\bigskip

In many interesting cases an answer can be given with the help of obstructions \ $\widetilde\omega (f_1, f_2)$ \ and \ $N (f_1, f_2)$ \ which we will now describe.

An analysis of the coincidence behaviour \ (of a suitable approximation) \ of \ $(f_1, f_2) : M \longrightarrow N \times N$ \ has led, in \cite{K 3}, to a normal bordism class
$$
\omega (f_1, f_2) \ = \ [C, g, \overline g] \ \in \ \Omega_{m -n} (M; \varphi)
\tag 1.2
$$
where \ $C$ \ is essentially the coincidence locus and \ $\overline g$ \ is a vector bundle isomorphism which describes the stable normal bundle of \ $C$ \ in terms of a pullback of the virtual coefficient bundle
$$
\varphi \ := \ f^*_1 (TN) - TM
\tag 1.3
$$
over \ $M$.  (Thus e.g.\ if \ $M$ \ and \ $N$ \ are stably parallelized, then  \ $\omega (f_1, f_2)$ lies in the framed bordism group \ $\Omega^{fr}_{m -n} (M)$).

In the present paper we will study a considerably sharper invariant based on the work of Hatcher and Quinn \cite{HQ}. Consider the commuting diagram
%% Einschub fuer Grafikeinbindung
$$
% eventuell Breite einstellen
% \EPSFxsize=5.5cm (andere Angabe)
\BoxedEPSF{diag1.eps}
\tag 1.4
$$
%% Einschub Ende
where \ $E (f_1, f_2)$ \ \ \footnote"*)"{ \ this differs slightly from the notation of Hatcher and Quinn.} consists of all pairs \ $(x, \theta)$ \ such that \ $x \in M$ \ and \ $\theta : [0, 1] \longrightarrow N$ \ is a continuous path with \ $\theta (0) = f_1 (x)$ \ and \ $\theta (1) = f_2 (x)$; \ $pr$ \ denotes the obvious projection and \ $\widetilde g$ \ is the natural lifting which adds the constant path at \ $f_1 (x) = f_2 (x)$ \ to \ $g (x) = x \in C$. Now, in general the topological space \ $E (f_1, f_2)$ \ is not pathconnected; in fact, its pathcomponents correspond bijectively to Reidemeister classes, i.e.\ to the classes in \ $\pi_1 (N)$ \ with respect to a certain relation depending on \ $f_1$ \ and \ $f_2$ \ (see prop. 2.1 below). Thus the normal bordism class
$$
\widetilde\omega (f_1, f_2) := [C, \widetilde g, \overline g] \ \in \ \Omega_{m -n} (E (f_1, f_2); \widetilde\varphi) = \bigoplus_{A \in \pi_0 (E (f_1, f_2))} \Omega_{m-n} (A; \widetilde\varphi | A)
\tag 1.5
$$
(where
$$
\widetilde\varphi \ := \ pr^* (\varphi) \ )
\tag 1.6
$$
can be decomposed into its contributions
$$
\widetilde\omega_A (f_1, f_2) = [ \ C_A := \widetilde g^{-1} (A), \ \widetilde g | C_A, \ \overline g | \ ] \ \in \ \Omega_{m -n} (A; \widetilde\varphi | A)
\tag 1.7
$$
to the various pathcomponents A.

\definition{Definition 1.8} \ A pathcomponent \ $A$ \ of \ $E (f_1, f_2)$ \ is called \ {\it essential} \ if the corresponding direct summand of \ $\widetilde \omega (f_1, f_2)$ is nontrivial.

The \ {\it Nielsen coincidence number} \ $N (f_1, f_2)$ \ is the number of essential pathcomponents \ $A \in \pi_0 (E (f_1, f_2))$.
\enddefinition

Since we assume \ $M$ \ to be compact, \ $N (f_1, f_2)$ \ is finite.

The Nielsen number is a simple (but crude) measure for the non-triviality of the much more delicate invariant \ $\widetilde \omega (f_1, f_2)$ \ which lies, a priori, in a group varying with \ $f_1$ \ and \ $f_2$. However, clearly \ $\widetilde\omega (f_1, f_2)$ \ vanishes if and only if \ $N (f_1, f_2)$ \ does.

\proclaim{Theorem 1.9} \
Given any pair of maps \ $f_1, f_2 : M \to N$ \ we have:
\roster
\item"(i)" \
$N (f_1, f_2)$ \ depends only on the homotopy classes of \ $f_1$ \ and \ $f_2$;
\item"(ii)" \
$ N (f_1, f_2) \ = \ N (f_2, f_1)$;
\item"(iii)" \
$N (f_1, f_2) \ \le \ MCC (f_1, f_2) \ < \ \infty$ \ \ $($cf.\ $1.1)$; \newline if \ $n \ne 2$, then also \ $MCC (f_1, f_2) \ \le \ \# \pi_0 (E (f_1, f_2))$  \ $($i.e.\ $MCC (f_1, f_2)$ \ is bounded by the cardinality of the Reidemeister set of \ $(f_1, f_2)$, cf.\ $2.1)$.
\endroster
\endproclaim
\medskip

In particular, it follows that \ $N (f_1, f_2)$ \ must vanish if \ $(f_1, f_2)$ \ is loose (cf.\ 1.1).
This leads us to ask the more specific

\subhead Question
\endsubhead
\ Is the Nielsen number \ $N (f_1, f_2)$ \ the only looseness obstruction for \ $(f_1, f_2)$? More generally: is \ $N (f_1, f_2)$ \ equal to \ $MCC (f_1, f_2)$?

\bigskip
As we will see below (cf.\ 1.10, 1.13, and 1.15) the answer is positive in at least three interesting settings, but not always (cf.\ 1.17).

\proclaim{Theorem 1.10} \ Assume \ $m < 2n - 2$.

Then for all maps \ $f_1, f_2 : M^m \to N^n$ \ we have \ $MCC (f_1, f_2) = N (f_1, f_2)$.

In particular, \ $(f_1, f_2)$ \ is loose if and only if \ $N (f_1, f_2) = 0$.
\endproclaim

This follows  from the work of Hatcher and Quinn \cite{HQ} (see \S\ 4 below). Their proof explains also the crucial role of the paths \ $\theta$ \ which occur in \ $E (f_1, f_2)$ \ and in a nulbordism of \ $(C, \widetilde g, \overline g)$ \ (cf.\ 1.4 and 1.5): they yield the necessary homotopies  \ $f_i \sim f'_i, \ i = 1,2$,  \ via a generalized Whitney trick construction.

\example{Example I:  \ ${\boldkey m \ \bold = \ \boldkey n}$ \ (the classical case)} \ Here \ $\widetilde\omega (f_1, f_2)$ \ lies in
$$
\Omega_0 (E (f_1, f_2); \widetilde\varphi) = \bigoplus \Sb A \in \pi_0 (E (f_1, f_2)) \\ \text{with}\ \widetilde\varphi | A \\ \text{oriented} \endSb \Bbb Z \ \ \ \ \ \oplus \bigoplus \Sb A \in \pi_0 (E (f_1, f_2)) \\ \text{with}\ \widetilde\varphi | A \\ \text{non orientable}\endSb \Bbb Z_2
$$
and counts the (generically transverse) coincidence points with or without signs. (The normal bordism approach makes it clear from the outset which of these two possibilities apply; in proposition 5.2 below we will give an explicit natural criterion in terms of fundamental groups). At least for orientable \ $M$ \ and \ $N$ \ the resulting contribution to a pathcomponent \ $A$ \ is just the index of the corresponding Reidemeister (or Nielsen) class (cf.\ e.g.\ \cite{BGZ}, 3.5 -- 3.6) and our definition of the Nielsen number agrees with the classical one. Moreover, a modified version of theorem 1.10 (cf.\ also proposition 4.7 below and \cite{Br}) yields the classical \lq\lq Wecken theorem\rq\rq\ for closed smooth manifolds of dimension at least 3 (cf.\ \cite{B}, p.\ 12).

The central idea of Nielsen fixed point theory -- to interpret the Lefschetz number (which corresponds to the \ $\omega$-invariant, cf.\ 1.2) as the sum of indices of the various Nielsen classes -- is expressed here by the $0$-dimensional normal bordism group of a single space. For \ $m - n > 0$ \ and arbitrary \ $(f_1, f_2)$ \ our approach seems to be even better suited to capture relevant geometric aspects: \ $\widetilde\omega_A (f_1, f_2)$ \ (cf.\ 1.7) will, in general, reflect much more than the oriented or unoriented bordism class of the underlying partial coincidence manifold \ $C_A$ \ or some derived (co)homology classes; often the full combined information contained in \ $C_A$, the map \ $\widetilde g | C_A$ \ and, in particular, the \lq\lq twisted framing\rq\rq\ \ $\overline g|$ \ turns out to be decisive (as illustrated e.g.\ by the examples in \cite{K 3}).
\endexample

\example{Example II: \ $\boldkey{f}_{1} \ \bold = \ \boldkey{f}{_2} \ \bold{=:} \ \boldkey{f}$ \ (selfcoincidence)} \
Here the projection \ $pr$ (cf.\ 1.4) allows a global section involving constant paths. Therefore \ $\widetilde\omega (f, f)$ \ (cf.\ 1.5) is precisely as strong as the obstruction \ $\omega (f) := \omega (f, f)$ \ (cf.\ 1.2) which was discussed in detail in \cite{K 3}. In particular, \ $N (f, f)$ \ equals \ $0$ or $1$ \ according as \ $\omega (f)$ \ vanishes or not (since \ $M$ \ is assumed to be
connected). Also clearly \ $MCC (f, f) = 1$ \ except when \ $(f, f)$ \ is loose (and hence \ $MCC (f, f) = 0)$. \hfill $\blacksquare$
\endexample

\medskip
An important special case of the \ $\widetilde\omega$--invariant is the \ {\it refined (normal bordism) degree}
$$
\widetilde{\deg} (f) := \widetilde\omega (f, *) \ \ \in \ \ \Omega_{m -n} (E (f, *); pr^* (f^* (TN) - TM))
\tag 1.11
$$
of a map \ $f : M \longrightarrow N$. \ Here \ $*$ \ denotes a constant map; its choice is not very significant since any path in \ $N$ \ from \ $*$ \ to \ $*'$, say, induces a fiber homotopy equivalence \ $E (f, *) \cong E (f, *')$ \ and an isomorphism of the corresponding normal bordism groups which takes \ $\widetilde\omega (f, *)$ \ to \ $\widetilde\omega (f, *')$. In particular, \ $\widetilde\omega (f, *)$ \ is compatible with the transitive action of \ $\pi_1 (N, *)$ \ on \ $\pi_0 (E (f, *))$. Therefore the pathcomponents  of \ $E (f, *)$ \ are either all essential (i.e. $N (f, *) = \# \pi_0 (E (f, *))\ = \ \# (\pi_1 (N) / f_* (\pi_1 (M)))$; in this case also \ $MCC (f, *) = N (f, *)$ \ if $n \ne 2)$ \ or all inessential (i.e. $N (f, *) = 0$; \ this holds e.g.\ if \ $\pi_1 (N) / f_* (\pi_1 (M)) $ \ is infinite).

The invariant \ $\widetilde{\deg} (f)$ \ sharpens the  (normal bordism) degree
$$
\deg (f) := \omega (f, *) \ \ \in \ \ \Omega_{m -n} (M; f^* (TN) - TM) \ 
\tag 1.12
$$
which was discussed extensively in \cite{K 3}.

Considerations of degrees (or \lq\lq roots\rq\rq) are often very useful in the general coincidence setting,  e.g.\ when \ $n$ \ is low (and when therefore theorem 1.10 offers little insight).

\example{Example III:  \ maps into the unit circle}
\endexample

\proclaim{Theorem  1.13} \ Assume \ $N = S^1$. Then the Nielsen number is characterized by the identity
$$
(f_{1 *} - f_{2 *}) (H_1(M; \Bbb Z)) = N (f_1, f_2) \cdot H_1 (S^1; \Bbb Z) \ \ ,
$$
and we have
$$
MCC (f_1, f_2) \ = \ N (f_1, f_2)
$$
$($cf.\ definition $1.1)$.

Moreover, the following conditions are equivalent (without any dimension restriction)
\roster
\item"(i)" \ $f_1$ \ and \ $f_2$ \ are homotopic;
\item"(ii)" \ $(f_1, f_2)$ \ is loose $($cf.\ definition $1.1)$;
\item"(iii)" \ the Nielsen number \ $N (f_1, f_2)$ \ vanishes;
\item"(iv)" \ $0 = \widetilde\omega (f_1, f_2) \in \Omega_{m -1} (E (f_1, f_2); - pr_* (TM))$;
\item"(v)" \ $0 = \omega (f_1, f_2) \in \Omega_{m -1} (M; - TM)$; \ \ \ and
\item"(vi)" \ $0 = \mu (\omega (f_1, f_2)) \in H_{m -1} (M; \widetilde{\Bbb Z_M}) {\underset PD\to\cong} H^1 (M; \Bbb Z)$
\endroster
(where \ $\mu$ \ and \ $PD$ \ denote the obvious Hurewicz and Poincar\'e duality homomorphisms).
\endproclaim

In this very special situation there is no need to do any calculations of normal bordism invariants: since \ $S^1$ \ is a Lie group and an Eilenberg--MacLane space, the Nielsen number \ $N (f_1, f_2)$ \ can be computed easily via degree theory (cf.\ 1.11 -- 1.12) and singular (co--)homology (and actually may assume all nonnegative integer values, for suitable choices of \ $M, f_1$ \ and \ $f_2$).

The detailed discussion of this example in  \S\ 6 will also include a simple and explicit description of the fiber homotopy type of \ $E (f_1, f_2)$ \ (which turns out to vary considerably,
depending on \ $f_{1*} - f_{2*} : H_1 (M; \Bbb Z) \to H_1 (S^1; \Bbb Z)$). \hfill $\blacksquare$

\example{Example IV: \ maps from spheres to spheres} \ Let \ $M = S^m, \ N = S^n$ \ with \ $n \ge 2$, \ and pick \ $y_0 \in S^n$. Then \ $S^n$ \  is 1-connected and hence the Nielsen number \ $N (f_1, f_2)$ \ equals 0 or 1 according as \ $\widetilde\omega (f_1, f_2)$ \ vanishes or not. Thus we really have to study the \ $\widetilde\omega$--invariant in detail here. After a canonical identification, all \ $\widetilde\omega$--invariants and the refined degree \ $\widetilde{\deg} (f)$ \ (compare 1.11) lie in the same group, namely the framed bordism group  \ $\Omega^{fr}_{m -n} (\wedge (S^n, y_0))$ \ of the loop space of \ $(S^n, y_0)$ \ which is independent of \ $f_1, f_2, f \dots$.
\endexample

\proclaim{Theorem 1.14} \ Given dimensions \ $m$ \ and \ $n \ge 2$, the refined $($normal bordism$)$ degree determines (and is determined by) a homomorphism which fits into the commuting  diagram
$$
\CD
[S^m, S^n]    @>{\widetilde\deg}>>  \Omega^{fr}_{m -n} (\wedge (S^n, y_0)) \\
@A{\cong}AA                               @V{\cong}V{h = \bigoplus_{k \ge 1} h_k}V \\
\pi_m (S^n) @>{\Gamma := \oplus E^\infty \scirc \gamma_k}>>   \bigoplus_{k \ge 1} \ \pi^S_{m -1 -k (n -1)} \ .
\endCD
$$
Here \ $E^\infty \scirc \gamma_k$ denotes the stabilized \ $k^{th}$ James-Hopf invariant homomorphism into the indicated stable homotopy groups of spheres, $k = 1, 2, \dots;  \ \ h_k$ is a geometric homomorphism which measures $(k -1)$--fold selfintersections of suitable immersions.

Given  maps \ $f_1, f_2, f$ \ from \ $S^m$ \ to \ $S^n$, \ $h$  yields a decomposition of \ $\widetilde\omega (f_1, f_2)$ \ and \ $ \widetilde\deg (f)$ \ into a sequence of components
$$
\widetilde\omega_k (f_1, f_2) \ := \ h_k (\widetilde\omega (f_1, f_2)), \ \ \widetilde\deg_k (f) := h_k (\widetilde\deg (f)) \ \ \in \ \ \pi^S_{m-1-k (n -1)} \ ,
$$
$k = 1, 2, \dots,$ \ starting with the (Pontryagin-Thom isomorphism evaluated on the nonrefined) invariants \ $\omega (f_1, f_2)$ \ and \ $\deg (f)$ = Freudenthal suspension \ $E^\infty (f)$ \ of \ $f$, resp. We have for all \ $k \ge 1$
$$
\widetilde\omega_k (f_1, f_2) \  = \  \widetilde{\deg}_k (f_1) - (-1)^{k (n -1)}  \widetilde\deg_k (f_2) \ ;
$$
moreover, if \ $n \not\equiv k \equiv 0 (2)$ \ or if \ $n \equiv 0 (2)$ \ and \ $k \equiv 3$ or $4 (4)$, then \ $2 \widetilde\omega_k (f_1, f_2) = 0$.

We have always
$$
MCC (f_1, f_2) = \left\{
\aligned 0 \ \ \text{if} \ f_1 \sim a \scirc f_2 \ , \\
1 \ \ \ \text{otherwise} \ \ \ ;
\endaligned \right.
$$
here \ $a$ \ denotes the antipodal map on \ $S^n$.
\endproclaim

\proclaim{Corollary 1.15} \ Assume that the total stabilized James-Hopf homomorphism \ $\Gamma$ \ (cf.\ {\rm 1.14}) is injective on \ $\pi_m (S^n)$.

Then for any two maps \ $f_1, f_2 : S^m \longrightarrow S^n$ \ we have \ $MCC (f_1, f_2) = N (f_1, f_2)$. In particular  the following conditions are equivalent:
\roster
\item"(i)" \ $(f_1, f_2)$ \ is loose \ ;
\item"(ii)" \ $N (f_1, f_2) = 0$ \ ; \ \ and
\item"(iii)" \ $f_1$ \ is homotopic to \ $a \scirc f_2$ \ (where \ $a : S^n \longrightarrow S^n$ \ denotes the antipodal map).
\endroster

\endproclaim

By Freudenthal's theorem, the assumption of this corollary is satisfied in the stable range \ $m < 2n - 1$ \ (where \ $\Gamma$ \ consists only of the stable suspension \ $E^\infty$) and hence is less restrictive than the dimension condition in theorem 1.10. With the added help of (higher) James-Hopf invariants we can answer our original question also in many nonstable dimension settings.

\proclaim{Corollary 1.16} \
If \ $m - n \le 3$ \ and \ $n \ge 1$, \ then \ $\Gamma$ \ is injective on \ $\pi_m (S^n)$ \ and therefore the conclusion of corollary {\rm 1.15} holds.
\endproclaim

However, the diagram in theorem 1.14 can also lead to opposite results.

\proclaim{Corollary 1.17} \
If \ $n \ne 1, 3, 7$ \ is odd and \ $m = 2n - 1$, or if e.g. $(m, n) = (8,4), (9,4), (9,3), (10,4), (16, 8), (17,8), (10 + n, n)$ for \ $3 \le n \le 11$, \ or $(24,6)$, \ then there exists a map \ $f : S^m \longrightarrow S^n$ \ such that the pair \ $(f, $ constant map$)$ \ is not loose although its Nielsen number vanishes.
\endproclaim

\subhead{Remark 1.18}\endsubhead \ In a future paper we will also study the minimum number
$$
M C (f_1, f_2) \ := \ \min \{ \# C (f'_1, f'_2) | f'_1 \sim f_1, f'_2 \sim f_2\}
$$
of coincidence {\it points} \ which plays a central role in classical fixed point and coincidence theory (where \ $m = n$; cf.\ e.g.\ \cite{B} and \cite{BGZ}). However, when \ $m > n$ \ then \ $MC (f_1, f_2)$ \ is very often infinite, and it seems more natural and illuminating to investigate \ $MCC (f_1, f_2)$ (cf.\ 1.1).
\medskip

\subhead{Remark 1.19}\endsubhead \
In classical fixed point theory it took 57 years to disprove the so-called Nielsen conjecture (which says that the minimum number \ $MC$, cf.\ 1.18, agrees with the Nielsen number; cf.\ \cite{B},p.\ 12--14, and \cite{Br}). The counterexamples to its higher-codimensional analogon provided by corollary 1.17 could be viewed as an indication that our generalized Nielsen number is too weak. In order to strengthen it, we may keep track e.g.\ of the fact that the map \ $g$ (cf.\ 1.4) is an {\it embedding}, and of the {\it nonstable} \ vector bundle isomorphism \ $\nu (C (f_1, f_2), M) \simeq f^*_1 (TN)|$ \ (compare 4.3). This procedure sharpens e.g.\ the degree to make it injective on homotopy classes of maps into a sphere (without adding much new insight). However, it seems that in general the resulting stronger invariants and bordism sets become quite unmanageable.
\medskip

\subhead{Convention}\endsubhead \ A framing of a smooth immersion or embedding is a (nonstable) trivialization of its normal bundle; a framing of a smooth manifold is a {\it stable} \ trivialization of its tangent bundle.

\vskip7mm
\specialhead 2.\ \ The space \ $E (f_1, f_2)$
\endspecialhead

Let \ $P (N)$ \ denote the space of all continuous paths \ $\theta : [0, 1] \to N$, endowed with the compact-open topology. Then the projection \ $pr : E (f_1, f_2) \to M$ \ (cf.\ 1.4) is just the pullback of the starting point/end point fibration \ $P (N) \to N  \times N$ \ by the map \ $(f_1, f_2) : M \to N \times N$.

In the next result we assume that there is a coincidence point \ $x_0 \in M$ \ and we put \ $y_0 := f_1 (x_0) = f_2 (x_0) \in N$ \ and \ $\theta_0$ := constant path at \ $y_0$. (If no such point exists, \ $(f_1, f_2)$ \ is loose and our original question is answered). We will identify the fiber of \ $pr$\ at\ $x_0$  \ with the {\it loop space} \ $\Lambda (N, y_0)$ \ of paths in \ $N$ \ starting and ending at \ $y_0$. Denote the fiber inclusion by \ $\incl$.

\proclaim{Proposition 2.1}
The sequence of group homomorphisms
$$
%\aligned
%\cdots \to \ssize{\pi_{k +1} (M, x_0)} @>{f_{1*} - f_{2*}}>> \ssize{\pi_{k +1} (N, y_0)} @>{\incl_*}>> \ssize{\pi_k %(E (f_1, f_2), (x_0,  \theta_0}&\ssize{))}  @>{pr_*}>> \ssize{\pi_k (M, x_0)} {\to \cdots} \\
% &\cdots \ssize{\longrightarrow \pi_1 (M, x_0)}
%\endaligned
% eventuell Breite einstellen
\EPSFxsize=12.5cm % (andere Angabe)
\BoxedEPSF{diag4.eps}
$$
is exact. Moreover the map
$$
\incl_* \ : \ \pi_1 (N,  y_0) = \pi_0 (\Lambda (N, y_0)) \longrightarrow \pi_0 (E (f_1, f_2))
$$
induces a bijection from the so called \ {\rm Reidemeister set}\
$$
\left. R (f_1, f_2, x_0) := \pi_1 (N, y_0)\ \right/ \ \text{\rm Reidemeister equivalence}
$$
onto the set of pathcomponents of \ $E (f_1, f_2)$. (We call \ $[\theta], [\theta'] \in \pi_1 (N, y_0)$ \ {\rm Reidemeister equivalent} \ if \ $[\theta'] = f_{1 *} (\mu)^{- 1} \cdot [\theta] \cdot f_{2 *} (\mu)$ \ for some \ $\mu \in \pi_1 (M, x_0))$.

Given a loop \ $\theta \in \Lambda (N, y_0)$ \ and the corresponding pathcomponent \ $A_{\theta}$ \ of \ $(x_0, \theta)$ \ in \ $E (f_1, f_2)$, we have 
$$
pr_* (\pi_1 (A_{\theta}; (x_0, \theta))) = \{ \mu \in \pi_1 (M, x_0) \vert f_{2*} (\mu) = [\theta]^{-1} \cdot f_{1*} (\mu) \cdot [\theta]\} \ .
$$
\endproclaim

\demo{Proof}
We are dealing here basically with the exact homotopy sequence of the Hurewicz fibration \ $pr$. The special form of the boundary homomorphism and the surjectivity of \ $\incl_*$ \ follow from the obvious
\enddemo

\remark{Fact {\rm 2.2}}  Given \ $(x, \theta) \in E (f_1, f_2)$, any path \ $\mu$ \ in \ $M$ \ from \ $x$ to $x_0$ lifts canonically to a path in \ $E (f_1, f_2)$ \ from \ $(x, \theta)$ \ to \ $(x_0, (f_1 \scirc \mu)^{-1} \ast \theta \ast (f_2 \scirc \mu))$ (here $\ast$ stands for successive travelling through paths, beginning with \ $(f_1 \scirc \mu)^{- 1})$.  \hfill $\blacksquare$
\endremark

\subhead{Remark 2.3}\endsubhead \ 
Given two coincidence points \ $x_1, x_2 \in M$ \ of \ $f_1, f_2$, their \ $\widetilde g$-values
$$
\widetilde g (x_i) = (x_i, \ \text{constant path at} \ f_1 (x_i) = f_2 (x_i)),
$$
$i = 1,2$, lie in the same path component of \ $E (f_1, f_2)$, precisely if the points \ $x_1$ and $x_2$ \ are {\it Nielsen equivalent}, \ i.e.\ they are joined by a path \ $\sigma$ \ in $M$ \ such that \ $f_1 \scirc \sigma$ \ and \ $f_2 \scirc \sigma$ \ are homotopic in \ $N$ \ leaving endpoints fixed.
\medskip

Now let \ $\varphi$ \ be a virtual vector bundle over \ $M$; put \ $\widetilde\varphi := pr^{*} (\varphi)$.

\proclaim{Proposition 2.4}
(see \cite{HQ}, 3.1). \ Assume that \ $M$ \ is \ $(k +1)$-connected for some integer $k$. Then a choice of an orientation of \ $\varphi$ \ at some point \ $x_0 \in M$ \ and of paths connecting \ $f_1 (x_0)$ \ and \ $f_2 (x_0)$ \ to some point \ $y_0$ \ in \ $N$ \ determine an isomorphism
$$
i \ : \Omega^{fr}_k (\Lambda (N, y_0)) \ \longrightarrow \Omega_k (E (f_1, f_2); \widetilde{\varphi}) \ \ .
$$
\endproclaim

\demo{Proof}
The choices induce isomorphisms
$$
\Omega^{fr}_k (\Lambda (N, y_0)) \ \cong \ \Omega^{fr}_k (pr^{- 1} (x_0)) \ \cong \ \Omega_k (pr^{- 1} (x_0); \widetilde\varphi|) \ 
$$
and so does the fiber inclusion; this follows via a cell-by-cell argument applied to (projected) maps into \ $M$.
\enddemo

\vskip7mm
\specialhead 3.\ \ Homotopies and vector bundle isomorphisms
\endspecialhead

We will need the following well known fact in the construction of the invariant \ $\widetilde\omega (f_1, f_2)$ \ and for establishing its homotopy invariance and symmetry properties.

\proclaim{Lemma 3.1}
Let \ $P$ \ be a $CW$-complex and \ $\xi$ \ an \ $n$-dimensional vector bundle over a space  \ $Q$. Then any homotopy \ $H : P \times I \ \longrightarrow \ Q$ between \ $h_i = H ( \ , i), \ i = 0,1$, determines a vector bundle isomorphism \ $h_0^* (\xi) \cong h^*_1 (\xi)$ \ over \ $P$, canonical up to regular homotopy.
\endproclaim

Indeed, the isomorphism bundle \ $\Iso (H^* (\xi), (h_0 \scirc \text{first projection})^* (\xi))$ \ over \ $P \times I$ \ (with fiber \ $GL (n, \Bbb R))$ \ has the required homotopy lifting property.

In particular, homotopies \ $F_i : M \times I \longrightarrow N$ \ from \ $f_i$ \ to \ $f_i', \ i = 1, 2, $ \ induce not only a fiber homotopy equivalence
$$
E (f_1, f_2) \ @>{ \ \sim \ }>> \ E (f'_1, f'_2)
\tag 3.2
$$
(which maps \ $(x, \theta)$ \ to \ $(x, (F_1 (x, -))^{- 1} \ast \theta \ast F_2 (x, -))$) \ but also a canonical isomorphism
$$
\Omega_* (E (f_1, f_2); pr^* (\varphi)) \ @>{ \ \cong \ }>> \ \Omega_* (E (f'_1, f'_2); pr^{'*} (\varphi'))
\tag 3.3
$$
where \ $\varphi = f^*_1 (TN) - TM$ \ (as in 1.3) and  \ $\varphi' := f'_{1*} (TN) - TM$. This isomorphism remains unchanged when we deform the homotopies \ $F_i$ \ while leaving them fixed on \ $M \times \{0, 1\}$.

Similarly the homeomorphism 
$$
E (f_1, f_2) \cong E (f_2, f_1), \ \ \ (x, \theta) \longrightarrow (x, \theta^{- 1})
$$
and the tautological homotopy between \ $f_i \circ pr, \ i = 1, 2$, determine canonical isomorphisms
$$
\aligned
\Omega_* (E (f_1, f_2); pr^* (\varphi)) &\cong \ \Omega_* (E (f_1, f_2); pr^* (f_2^* (TN) - TM)) \\
                                         &\cong \ \Omega_* (E (f_2, f_1); pr^* (f_2^* (TN) - TM)).
\endaligned
\tag 3.4
$$
\subhead{Remark 3.5}\endsubhead
Such symmetries are due to the paths occuring in \ $E (f_1, f_2)$ \ and to the resulting homotopy. In general, they are not available in the corresponding normal bordism groups of \ $M$. E.g. \ if \ $TM$ \ and \ $f^*_1 (TN)$ \ are orientable, but \ $f^*_2 (TN)$ \ is not, then
$$
\Omega_0 (M; f_1^* (TN) - TM) \cong \Bbb Z \ \not\cong \ \Bbb Z_2 \cong \Omega_0 (M; f^*_2 (TN) - TM) \ .
$$

\vskip7mm
\specialhead 4.\ \ The invariants \ $\widetilde\omega (f_1, f_2), \ N (f_1, f_2)$ \ and \ $\omega (f_1, f_2)$.
\endspecialhead

In this section we construct and discuss our coincidence invariants. Moreover, we describe a technique to make Nielsen classes pathconnected. In particular, we prove theorems 1.9 and 1.10.  We also sketch the definition of certain secondary obstructions in normal bordism.

If the map \ $(f_1, f_2) : M \longrightarrow N \times N$ \  is smooth and transverse to the diagonal
$$
\Delta \ := \ \{ (y, y) \  \in \  N \times N \ | \ y \in N \}
$$
then the coincidence locus
$$
C \ = \ C (f_1, f_2) \ := \ \{ x \in M \ | \ f_1 (x) = f_2 (x) \} \ =  (f_1, f_2)^{- 1} (\Delta)
\tag 4.1
$$
is a closed smooth \ $(m - n)$-dimensional manifold, equipped with the map
$$
\widetilde g \ : \ C \ \longrightarrow \ E (f_1, f_2)
\tag 4.2
$$
which sends \ $x \in C$ \ to \ $(x,$ constant path at $f_1 (x) = f_2 (x))$, and with a stable tangent bundle isomorphism
$$
\overline g \ : \ TC \oplus \widetilde g^* (pr^* (f_1^* (TN))) \cong \widetilde g^* (pr^* (TM))
\tag 4.3
$$
(since the normal bundle \ $\nu (\Delta, N \times N)$ \ of \ $\Delta$ \ in \ $N \times N$ \ is canonically isomorphic to the pullback of the tangent bundle \ $TN$ \ under the first projection \ $p_1$).

If \ $f_1$ and $f_2$ \ are arbitrary continuous maps, apply the preceding construction to a smooth map \ $(f'_1, f'_2)$ \ which approximates \ $(f_1, f_2)$ \ and is transverse to \ $\Delta$. Then there is a canonical isomorphism as in 3.3 induced by any sufficiently small homotopy from \ $(f_1, f_2)$ \ to \ $(f'_1, f'_2)$.

In any case the resulting triple \ $(C, \widetilde g,  \overline g)$ \ determines a well-defined normal bordism class
$$
\widetilde\omega (f_1, f_2) \ \in \ \Omega_{m -n} (E (f_1, f_2); \ \widetilde\varphi) \ .
\tag 4.4
$$

Arbitrary (possibly \lq\lq large\rq\rq) homotopies of \ $f_i$ \ yield a homotopy equivalence as in 3.2 and, in particular, a bijection of essential path components (cf.\ 1.8); indeed, the isomorphism 3.3 is compatible with \ $\widetilde\omega$. This proves the homotopy invariance of the Nielsen number claimed in 1.9 (i).

Similarly, \ $\widetilde\omega$ \ is compatible with the isomorphism 3.4 which, however, must first be composed with a suitable involution; this is needed since over \ $N \cong \Delta$ \ the composite of the obvious isomorphisms
$$
TN \ \cong \ p_1^* (TN) \ \cong \ \nu (\Delta, N \times N) \ \cong \ p^*_2 (TN) \ \cong \ TN
\tag 4.5
$$
equals \ $- \id_{TN}$. The second claim in theorem 1.9 follows. 
\medskip

It remains to estimate the number of pathcomponents of the coincidence locus \ $C (f_1, f_2)$ \ (cf.\ 4.1) when \ $f_1, f_2$ are arbitrary maps (without smoothness or transversality assumptions). Embed \ $N$ \ as a closed submanifold of \ $\Bbb R^{2 n+1}$ \ and compose a linear homotopy in \ $\Bbb R^{2 n+1}$ \ with a tubular neighbourhood projection onto \ $N$ \ to obtain a map
$$
\tau \ : \ U \times I \ \longrightarrow \ N \ \ \ \text{such that} \ \ \tau (y_0, y_1, i) \ = \ y_i, \ \ i = 0, 1 \ ,
$$
for all \ $(y_0, y_1)$ \ in a suitable neighbourhood \ $U$ \ of the diagonal \ $\Delta$ \ in \ $N \times N$. This allows us to extend the lifting \ $\widetilde g$ \ (cf.\ 1.4) to the map
$$
\widetilde G \ \ \  : \ \ \ V  := (f_1, f_2)^{-1} (U) \ @>{ \ \ \ \ \ }>> \ E (f_1, f_2) \ ,
\tag 4.6
$$
$\widetilde G (x) := (x, \tau (f_1 (x), f_2 (x), - )$, defined in the neighbourhood \ $V$ \ of the coincidence locus \ $C (f_1, f_2) = (f_1, f_2)^{- 1} (\Delta)$ \ in \ $M$. The decomposition of \ $E (f_1, f_2) = \cup \ A$ into (open!) pathcomponents yields a corresponding decomposition \ $V = \cup \ V_A$ \ which is compatible with any small (smooth, transverse) approximation \ $(f'_1, f'_2)$ \ of \ $(f_1, f_2)$. If \ $V_A \cap \ C (f_1, f_2) = \varnothing$, we may assume that the corresponding part of \ $C (f'_1, f'_2)$ \ is also empty. Even if \ $V_A \cap \ C (f_1, f_2) \ne \varnothing$, the corresponding part of \ $C (f'_1, f'_2)$ \ may be nonessential. Thus clearly
$$
\# \pi_0 (C (f_1, f_2)) \ \ge \ N (f'_1, f'_2) \ = \ N (f_1, f_2) \ .
$$
Moreover, the compact manifold \ $C (f'_1, f'_2)$ \ has only finitely many pathcomponents.

The last claim in theorem 1.9 follows from

\proclaim{Proposition 4.7} \ Assume \ $n \ne 2$. Given \ $(f_1, f_2)$, there exists a pair of maps \ $(f'_1, f'_2)$ \ and a homotopy \ $(f_1, f_2) \sim (f'_1, f'_2)$ \ relating each pathcomponent \ $A \in \pi_0 (E (f_1, f_2))$ \ $($with Nielsen class \ $C_A \subset C (f_1, f_2)$, cf.\ {\rm 1.7, 2.3,}\ and {\rm 3.2)} to a corresponding pathcomponent \ $A' \in \pi_0 (E (f'_1, f'_2))$ $($with Nielsen class \ $C_{A'} \subset C (f'_1, f'_2))$, in such a way that
\roster
\item"(i)" \ each Nielsen class \ $C_{A'} \subset C (f'_1, f'_2)$ \ is pathconnected $($or empty$)$;
\item"(ii)" \ if \ $C_A$ \ is finite, then \ $C_{A'}$ \ consists of at most a single point; \ and
\item"(iii)" \ if \ $C_A$ \ is empty, then so is \ $C_{A'}$.
\endroster
\endproclaim

\demo{Proof} \
If \ $m < n$ \ or \ $N = \Bbb R$, then \ $(f_1, f_2)$ \ is loose, and after a suitable homotopy every \ $C_{A'}$ \ is empty. The case \ $N = S^1$ \ will be discussed very explicitly in \ \S\ 6 (at the end of the proof of theorem 1.13; also note: when \ $m \ge 2$ \ every {\it isolated} \ coincidence point \ $x$ \ can be removed here by a deformation with small support near $x$. So we may assume that \ $m \ge n \ge 3$.

We will use the technique of the previous proof (cf.\ the lines following 4.6) to ensure that sufficiently small deformations do not cause new nonempty Nielsen classes. Thus, given \ $x \in C (f_1, f_2)$, we may assume that \ $(f_1, f_2)$ \ is smooth and transverse to the diagonal \ $\Delta \subset N \times N$ \ near \ $x$; \ or, if the coincidence point \ $x$ \ is isolated, then (in terms of local coordinates) \ $f_1 - f_2$ \ is linear on each ray starting from \ $x$. In either case in the end \ $C (f_1, f_2)$ \ and hence each Nielsen class \ $C_A$ \ consists only of finitely many pathcomponents. We decrease their number by iterating the following procedure.

If the points \ $x_0, x_1 \in C_A$ \ lie in different pathcomponents join them by a path \ $\sigma : I \to M$ \ which otherwise avoids \ $C (f_1, f_2)$ \ and such that \ $f_1  \scirc \sigma \sim f_2 \scirc \sigma \ \ \rel \{ 0, 1\}$ \ via a homotopy \ $F : I \times I \to N$. Since \ $m \ge n \ge 3$ \ we may assume that \ $\sigma$ \ and \ $f_1 \scirc \sigma$ \ are smooth embeddings and that the transverse intersection of \ $F ((\varepsilon, 1-\varepsilon) \times (0,1])$ \ with \ $f_1 \scirc \sigma (I)$ \ consists of at most finitely many points of the form \ $F (t, s) = f_1 \scirc \sigma (t')$ \ such that \ $t \ne t'$. Now deform \ $f_1, f_2$ \ simultaneously near a small tubular neighbourhood \ $I \times B^{n -1}$ \ of \ $\sigma (I) = I \times \{ 0\}$ \ until \ $f_i (t, x') = f_i \scirc \sigma (t), \ \ i = 1,2$, \ for all \ $t \in I$ \ and \ $x'$ \ in the unit ball \ $B^{n -1}$; then replace \ $f_2 (t, x')$ by \ $F (t, \Vert x'\Vert)$. This modification creates new coincidence points only along the whole arc \ $\sigma (I)$ \ and, in addition, possibly near its endpoints. We complete the proof of proposition 4.7 by applying the following result to suitable compact balls containing the arcs \ $f_1 \scirc \sigma (I) = f_2 \scirc \sigma (I)$ \ and \ $\sigma (I)$.
\enddemo

\proclaim{Lemma 4.8} \
Given arbitrary continuous maps \ $f_1, f_2 : M \to N$, a compact $n$-ball \ $Q \subset N$ \ with collar, and a compact subset \ $P \subset f^{-1}_1 (Q) \cap f^{- 1}_2 (Q) \subset M$, there exist homotopic maps \ $f'_1 \sim f_1$ \ and \ $f'_2 \sim f_2$ \ with coincidence locus
$$
C (f'_1, f'_2) \ \ = \ \ C (f_1, f_2) \ \ \cup \ \ P \ \ .
$$
If, in addition, \ $P$ \ is a compact \ $m$-ball in \ $M$ \ with collar, then after a further homotopy \ $C (f'_1, f'_2)$ \ can be made to be homeomorphic to the quotient space \ $(C (f_1, f_2) \cup P) / P$.

$($No restrictions are imposed here on \ $m = \dim M$ \ and \ $n = \dim N)$.
\endproclaim

\demo{Proof} \
The assumption on \ $Q$ \ means that \ $Q$ \ is the image of the compact \ $n$-ball \ $B^n (1) \subset \Bbb R^n$ \ of  radius $1$ \ under a homeomorphism \ $B^n (1 + \varepsilon) \cong \widehat Q \subset N$ \ for some \ $\varepsilon > 0$. Let $* \in Q$ \ be the centre point and let \ Cone $(Q) \subset N \times I$ \ denote the cone with basis \ $Q \times \{ 0\}$ \ and top vertex \ $z_0 := (*, 1)$. By shrinking \ $Q$ \ and stretching the collar radially we can construct a continuous family \ $(h_s), \ s \in [0,1)$, of selfhomeomorphisms of \ $N$ \ which starts with \ $h_0 = \id$ \ and converges to a selfmap \ $h_1$ \ of \ $N$ \ taking \ $N - Q$ \ homeomorphically to \ $N - \{*\}$ \ and such that \ $h_1 (Q) = \{ *\}$.

Next choose a continuous function \ $\delta : M \to [0, 1]$ \ satisfying \ $\delta^{-1} (\{ 1\} ) = P$. The required maps \ $f'_i$ \ can be defined by
$$
f'_i (x) \ = \ h_{\delta (x)} (f_i (x))  \  , \ \ \  x\in M, \ \ \ i = 1, 2 \ .
$$ 

If \ $P \subset M$ \ is also an \ $m$-ball with collar, construct a similar homotopy \ $(\overline h_s)$ \ of selfmaps of \ $M$. \ Then the maps \ $f'_i$ \ factor through \ $\overline h_1, \ \ f'_i = f''_i \scirc \overline h_1 (\sim f''_i \scirc \overline h_0 = f''_i, \ i = 1, 2,$ \ and \ $C (f''_1, f''_2) = (C (f_1, f_2) \cup P) / P$. \hfill $\blacksquare$
\enddemo
\medskip

Clearly the invariants \ $\omega (f_1, f_2)$ \ (defined in \cite{K 3}; see also 1.2) and \ $\widetilde\omega (f_1, f_2)$ \ (cf.\ 4.4) are related by the equation
$$
\omega (f_1, f_2) \ = \ pr_* (\widetilde\omega (f_1, f_2))
\tag 4.9
$$
where the homomorphism
$$
pr_* \ : \ \Omega_* (E (f_1, f_2); \widetilde\varphi) \ \longrightarrow \ \Omega_* (M; \varphi)
\tag 4.10
$$
is induced by the fiber projection \ $pr$ \ (cf.\ 1.4, 1.3, and 1.6). In general it leads to a loss not only of information but also of symmetry: while \ $\widetilde\omega (f_1, f_2)$ \ and \ $\widetilde\omega (f_2, f_1)$ \ are always equally strong, the invariants \ $\omega (f_1, f_2)$ \ and \ $\omega (f_2, f_1)$ \ need not be (compare remark 3.5).

In the very special case \ $f_1 = f_2 =: f$ \ of example II \ $pr$ has a global section \ $s$ \ such that \ $\widetilde\omega (f, f)$ \ equals \ $s_* (\omega (f, f))$ \ and hence is precisely as strong as the selfcoincidence invariant \ $\omega (f)$ \ studied in \cite{Ko 3}.

In general, however, it may happen that \ $\omega (f_1, f_2) = 0$, but \ $\widetilde\omega (f_1, f_2)$ \ gives rise to interesting secondary obstructions of the form \ $h (\widetilde\omega (f_1, f_2))$, where \ $h$ \  is a well-defined homomorphism on the kernel of \ $pr_*$. E.g.\ if \ $M$ and $N$ \ are stably parallelizable (and hence \ $\ker pr_*$ \ lies in the framed bordism group \ $\Omega^{fr}_* (E (f_1, f_2)))$, we may construct a homomorphism
$$
h_2 \ : \ \ker pr_*  \ \longrightarrow \ \widetilde\Omega_{*+1}^{fr} (N) / (f_{1*} - f_{2*}) (\Omega^{fr}_{* +1} (M))
\tag 4.11
$$
in the spirit of 2.1 as follows. Given \ $c \in \ker pr_*$, pick a framed singular manifold \ $C @>{\widetilde g = (g, \theta)}>> E (f_1, f_2)$ \ representing \ $c$ \ and a  framed nulbordism \ $G : B \to M$ \ 
of \ $(C, g)$;  define \ $h_2 (c)$ \ to be the class of the closed framed manifold
$$
- B \ \cup_{\partial B = C \times \{0\}} \ C \times I \ \cup_{C \times \{1\} = \partial B} \ B
$$
together with the map \ $f_1 \scirc G \cup \theta \cup f_2 \scirc G$ \ into $N$. If in addition \ $N = S^n, n \ge 2$, then the secondary obstruction \ $h_2 (\widetilde\omega (f_1, f_2))$ \ lies in a quotient of the reduced bordism group \ $\widetilde\Omega^{fr}_{m - n+1} (S^n) \cong \Omega^{fr}_{m -2n+1} \cong \pi^S_{m -1 -2 (n -1)}$ \ and generalizes the invariant \ $\widetilde\omega_2 (f_1, f_2)$ \ discussed in theorem 1.14 and section 8.  \hfill $\blacksquare$
\medskip

Next we turn to the

\demo{Proof of theorem {\rm 1.10}} \
By transversality we may assume that \ $m \ge n$ \ and hence \ $n > 2$. In view of proposition 4.7 we only have to find a homotopy which relates each inessential Nielsen class \ $C_A$ \ to a corresponding \ {\it empty} \ class \ $C_{A'}$. For this purpose we apply the fibered theory of Hatcher and Quinn (see \S\ 4 of \cite{HQ}).

After an approximation \ $(f_1, f_2)$ \ is smooth and transverse to \ $\Delta$. Also we have the smooth, transverse, \lq\lq fiberwise\rq\rq embeddings
$$
M \ @>{i_Q := (\id, f_1)}>> \ M \times N \ @<{i_P := (\id, f_2)}<< \ M
$$
which preserve the obvious projections onto \ $M$. Moreover, the diagonal map \ $M \to M \times M$ \ induces a homomorphism of normal bordism groups which takes \ $\widetilde\omega (f_1, f_2)$ \ to the invariant \ $[i_P \pitchfork i_Q]$ \ of Hatcher and Quinn.

Now replace each part \ $C_A$ \ of the coincidence locus \ $C (f_1, f_2)$ \ (cf. 4.1) which corresponds to an inessential path component of \ $E (f_1, f_2)$ \ (i.e.\ the triple \ $(C_A, \widetilde g | C_A, \overline g|)$ \ is  nulbordant , cf.\ 1.7) by the empty set. According to theorem 4.2 in \cite{HQ} the resulting modified coincidence data can be realized as the coincidence data of \ $i_P$ and an embedding \ $i'_Q = \id \times f'_1$ \ which is fiber isotopic to \ $i_Q$. In particular, \ $f'_1$ is homotopic to \ $f_1$. \hfill $\blacksquare$
\enddemo

The dimension restriction in theorem 1.10 turns out to be unnecessary in a number of cases. For a better understanding  of one such case in \S\ 6 it is helpful to study precompositions.

\proclaim{Proposition 4.12} \
Let \ $L$ \ be a smooth closed connected manifold and \ $e : L \to M$ \ a map. Assume that in the diagram
$$
\pi_1 (L, \ell_0) \ @>{ \ e_* \ }>> \ \pi_1 (M, x_0 := e (\ell_0)) \ @>{ \ (f_1, f_2)_* \ }>> \ \pi_1 (N \times N, (f_1, f_2) (x_0)) 
$$
$(f_1, f_2)_* \scirc e_*$ \ and \ $(f_1, f_2)_*$ \ have equal images.
\endproclaim

Then \ $N (f_1 \scirc e, f_2 \scirc e) \ \le \ N (f_1, f_2)$.

\demo{Proof}
In view of 1.9 (i) we need only consider the case where \ $(f_1, f_2)$ \ and \ $e$ \ are smooth and transverse to the diagonal \ $\Delta$ \ and to the coincidence locus \ $C (f_1, f_2) = (f_1, f_2)^{- 1} (\Delta)$, resp. (cf.\ 4.1), and where in addition \ $N (f_1 \scirc e, f_2 \scirc e) > 0$ \ and \ $x_0 \in C (f_1, f_2)$. Then it follows from our assumption and from 2.1 that the map
$$
E (f_1 \scirc e, f_2 \scirc e) \ \longrightarrow \ E (f_1, f_2)
$$
determined by \ $e$ \ makes each path component \ $A$ of $E (f_1, f_2)$ \ correspond to a unique path component \ $A'$ \ of \ $E (f_1 \scirc e, f_2 \scirc e)$. Now, if \ $A$ \ is inessential then so is \ $A'$. Indeed, let \ $C_A$ \ and \ $C_{A'}$ \ denote the corresponding parts of the coincidence manifolds \ $C (f_1, f_2)$ \ and \ $C (f_1 \scirc e, f_2 \scirc e) = e^{- 1} (C (f_1, f_2))$. Then a nulbordism of \ $C_A$ \ yields a nulbordism of \ $C_{A'}$ \ by transverse intersection with \ $e$ \ in \ $M$. \hfill $\blacksquare$
\enddemo

For later use we note also one of the many interrelations with degree theory (compare 1.11).

\proclaim{Proposition 4.13} \ 
If \ $N$ \ is a Lie group (and hence there is a multiplication of maps into \ $N$), then
$$
N (f_1, f_2) \ = \ N (f_1  \cdot f_2^{-1}, \ \text{constant map} \ ) \ 
$$
and this Nielsen number equals \ $0$ \ or \ $\# (\pi_1 (N) / (f_{1*} - f_{2*}) (\pi_1 (M)))$. Moreover,
$$
MCC (f_1, f_2) \ = \ MCC (f_1 \cdot f_2^{-1}, \ \text{constant map} \ )  \ \ .
$$
\endproclaim

\demo{Proof} \ 
The diffeomorphism
$$
(N \times N, \Delta) \ \longrightarrow \ (N \times N, 1 \times N) \ ,
$$
$(y_1, y_2) \longrightarrow (y_1 \cdot y_2^{- 1}, y_2)$, and the natural homeomorphism \ $E (f_1, f_2) \cong E (f_1 \cdot f_2^{- 1}, 1)$ \ lead to an isomorphism of normal bordism groups which takes \ $\widetilde\omega (f_1, f_2)$ \ to \ $\widetilde\omega (f_1 \cdot f_2^{- 1}, 1) = \widetilde\deg (f_1 \cdot f_2^{-1})$. In particular, the corresponding Nielsen numbers agree. Also the transitive action of \ $\pi_1 (N)$ \ on \ $\pi_0 (E (f_1 \cdot f_2^{- 1}, 1)) (\approx \coker (f_1 \cdot f_2^{-1})_*$, cf.\ 2.1) is compatible with \ $\widetilde\deg (f_1 \cdot f_2^{- 1})$. Hence either all pathcomponents of \ $E (f_1, f_2)$ \ are essential or all are inessential.

If the pairs of maps \ $(f_1 \cdot f_2^{- 1}, 1)$ \ and \ $(g_1, g_2)$ \ are homotopic, then so are \ $(f_1, f_2)$ \ and \ $(f_1 \cdot f_2, g_2 \cdot f_2)$. Since \ $(f_1, g_2)$ \ and \ $(g_1 \cdot f_2, g_2 \cdot f_2)$ \ have the same coincidence set we conclude that \ $MCC (f_1, f_2) \le MCC (f_1 \cdot f_2^{ - 1}, 1)$. The full equality follows similarly.
\enddemo

\vskip7mm
\specialhead 5.\ \ Classical Nielsen coincidence theory
\endspecialhead

In this section we consider the case when \ $M$ \ and \ $N$ \ have the same dimension. Then
$$
\widetilde\omega (f_1, f_2) \ = \ \{ \widetilde\omega_A (f_1, f_2)\} \  \ \in \ \ \oplus_{A \in \pi_0 (E (f_1, f_2))} \ \Omega_0 (A; \widetilde\varphi | A)
$$
(cf.\ 1.5 and 1.7).

Now, a \ $0$-dimensional normal bordism group of any pathconnected (and locally pathconnected) space  is isomorphic to \ $\Bbb Z$ \ or \ $\Bbb Z_2$ \ according as the coefficient bundle is orientable or not, and this is usually expressed by the first Stiefel-Whitney class. In the situation at hand, \ $\widetilde\varphi | A = (pr | A)^* (\varphi)$ \ is orientable if and only if
$$
w_1 (\varphi) \ = \ w_1 (M) + f_1^* (w_1 (N)) \ \in \ H^1 (M; \Bbb Z_2)
\tag 5.1
$$
(cf.\ 1.3 and 1.6) vanishes when evaluated on the image of the fundamental group of \ $A$ \ under \ $pr_*$. In view of proposition 2.1 we obtain

\proclaim{Proposition 5.2} \
Let \ $x_0 \in M$ \ be a coincidence point with \ $y_{0} :=  f_1 (x_0) = f_2 (x_0)$. Given a pathcomponent \ $A$ \ of \ $E (f_1, f_2)$, pick \ $[\theta] \ in \ \pi_1 (N, y_0)$ \ such that \ $(x_0, \theta) \in A$ \ (i.e. \ $A = A_{\theta}$ \ corresponds to the class of \ $[\theta]$ \ in the Reidemeister set).

Then \ $\Omega_0 (A; \widetilde\varphi | A) \cong \Bbb Z$ \ if and only if
$$
w_1 (M) (\mu) \ = \ f^*_1 (w_1 (N)) (\mu)
$$
for all \ $\mu \in \pi_1 (M, x_0)$ \ such that \ $f_{2*} (\mu) = [\theta]^{-1} \cdot f_{1*} (\mu) \cdot [\theta]$. Otherwise \ $\Omega_0 (A; \widetilde\varphi | A)$ $  \cong  \Bbb Z_2$.
\endproclaim

If \ $\pi_1 (N, y_0)$ \ is abelian, then this criterion is independent of the pathcomponent \ $A$ \ and we have
$$
\Omega_0 (E (f_1, f_2); \widetilde\varphi) \cong \left\{
\aligned \oplus_A & \Bbb Z \ \ \ \ \text{if} \ w_1 (\varphi) \equiv 0 \ \ \text{on} \ \ \ker (f_{2 *} - f_{1 *}) \ ; \\
\oplus_A & \Bbb Z_2 \ \ \ \text{else} \ \ \ .
\endaligned \right.
$$
However, if \ $\pi_1 (N, y_0)$ \ is not commutative, then both \ $\Bbb Z$ \ and \ $\Bbb Z_2$ \ may occur as direct summands (possibly corresponding to inessential Nielsen classes), e.g.\ in some cases when \ $M$ \ is the Klein bottle and \ $N$ \ is the punctured torus.
\hfill $\blacksquare$

\medskip
Now consider the generic situation where \ $(f_1, f_2)$ is smooth and transverse to the diagonal \ $\Delta$ \ in \ $N \times N$. Then  a coincidence point \ $x$ \ can contribute a well-defined \ {\it integer} \ index $\pm 1$ \ (the intersection number of \ $(f_1, f_2)$ \ with \ $\Delta$ \ at \ $x)$ \ to \ $\widetilde\omega (f_1, f_2)$  \ if and only if \ $w_1 (\varphi)$ \ vanishes on all \ $\mu \in \pi_1 (M, x)$ \ such that \ $f_{1 *} (\mu) = f_{2 *} (\mu)$ \ (the sign depends then on the choice of an orientation of \ $\widetilde\varphi$ \ on the path component of \ $(x$, constant path at \ $f_1 (x) = f_2 (x))$. E.g.\ in fixed point theory (where \ $f_1 = \id$ \ and hence \ $\varphi = 0$), or, more generally, if \ $f_1$ \ is \ {\it orientation true} \ (cf.\ \cite{BGZ}, p.\ 49), i.e.\ if we can (and do) choose an orientation already of \ $\varphi$ \ (not just of \ $\widetilde\varphi|A$), we can attach well-defined \ {\it integer} \  indices to all coincidence points.

Whether indices lie in \ $\Bbb Z$ \ or in \ $\Bbb Z_2$, \ $\widetilde\omega (f_1, f_2)$ \ assigns to every path component, i.e.\ to every Nielsen class (cf. 2.3) of coincidence points, its (total) index in a very natural way,  and our definition of Nielsen numbers coincides with the classical one (e.g.\ as in \cite{BGZ}, 3.5 and 3.6) at least if \ $M$ and $N$ \ are orientable.

In case \ $M = N$ \ has dimension at least \ $3$ \ and \ $f$ \ is a selfmap of \ $M$ \ we can easily modify the proof of theorem 1.10 (using proposition 4.7) to show that \ $N (f, \id)$ \ agrees with the minimum number
$$
MC (f, \id) \ := \  \min \{ \# C (f'_1, f'_2) | f'_1 \sim f, f'_2 \sim \id \}
$$
of coincidence {\it points} \ (not just pathcomponents) or, by \cite{BR}, with the minimum number
$$
MF (f) \ := \ \min \{ \# C (f', \id) | f' \sim f \}
$$
of fixed points. This is the classical Wecken theorem for closed smooth manifolds (cf.\ \cite{B}, p.\ 12).

This section used the language of $0$-dimensional bordism theory. Techniques for calculating normal bordism groups also in dimensions \ $m - n = 1, 2, 3$ \ can be found e.g.\ in \cite{K 1}, \S\ 9.

\vskip7mm
\specialhead 6.\ \ The case when $N$ is 1-dimensional
\endspecialhead

Because of the dimension restriction \ $m < 2n - 2$ \ theorem 1.10 is uninteresting for low \ $n$. In this section we will study the case \ $n = 1$ \ by different methods and answer our original question explicitly.

If \ $N = \Bbb R$ \ there is no problem: obviously \ $f_1$ and $f_2$ \ are homotopic to different (and hence coincidence free) constant maps.

For the remainder of this section let \ $N$ \ be the circle \ $S^1$. Now \ $S^1$ \ enjoys three very special properties. First of all, it is a Lie group and we can use proposition 4.13. As a consequence we need to prove theorem 1.13 only for pairs of the form \ $(f, *)$ \ where \ $*$ \ stands for a constant map with value, say, \ $y_0 = f (x_0)$ \ for some \ $x_0 \in M$. Let \ $d$ \ be the nonnegative integer characterized by the equation 
$$
f_* (\pi_1 (M, x_0)) \ = \ d \cdot \pi_1 (S^1, y_0) \ \cong \ d \Bbb Z .
$$
Then, on one hand, we have
$$
N (f, *) \ \le \ d
$$
since we can index \ {\it all} \ the  path components of \ $E ( f, *)$ \ (essential or not) by the Reidemeister set \ $R (f, *, x_0) = \coker f_*  \cong \Bbb Z/d \Bbb Z$ \ (cf. prop. 2.1). On the other hand, if we pick a loop \ $e : (S^1, 1) \to (M, x_0)$ \ such that \ $f \circ e : S^1 \to S^1$ \ has degree \ $d$, we obtain
$$
d \ = \ N (f \circ e, *) \ \le \ N (f, *) \ .
$$
Indeed, the inequality to the right follows from proposition 4.12; moreover, \ $f \scirc e$ \ is homotopic to the standard selfmap of the circle of degree \ $d$, whose \ $d$ \ roots clearly lie in pairwise distinct Nielsen classes (compare 2.3). This yields one possible proof of the first claim in theorem 1.13.

An alternative proof and the third claim follow from another special property of the circle: \ $S^1$ \ is an Eilenberg--MacLane--space for \ $H^1 ( \ ; \Bbb Z)$. Thus the homotopy class of \ $f$ \ is already determined by \ $f_* : H_1 (M; \Bbb Z) \to H_1 (S^1; \Bbb Z)$. Also we obtain the indicated isomorphisms in the commuting diagram
%% Einschub fuer Grafikeinbindung
$$
% eventuell Breite einstellen
\EPSFxsize=7cm % (andere Angabe)
\BoxedEPSF{diaghm.eps}
$$
%% Einschub Ende

%$$
%\CD
%[M, S^1]    @>{ \ \cong \ }>>     H^1(M, \Bbb Z) \\
%@VV{\deg}V                          @VVV \\
%\Omega_{m - 1} (M; -TM)  @>>{\ \mu \ }>  H_{m -1} (M; \widetilde {\Bbb Z}_M) \ .
%\endCD
%$$

Finally recall that \ $S^1$ \ is also a Thom space. Thus we can describe a deformation of a generic map \ $f : M \to S^1$ \ by exhibiting a bordism of the cooriented $1$-codimensional submanifold \ $C = f^{- 1} \{*\}$ \ of \ $M$.

Assume that \ $m \ge 2$ \ and that there are different pathcomponents \ $C_0$ \ and \ $C_1$ \ of \ $C$ \ belonging to the same Nielsen class (cf.\ 2.3). They can be joined by a smoothly immersed curve \ $\sigma : [ 0, 1] \to M$ \ which is transverse to \ $C$ \ and such that \ $f \scirc \sigma$ \ is nulhomotopic \ $\rel \{ 0, 1\}$. Let the crossings of \ $\sigma$ \ with \ $C$ \ have opposite signs at \ $\sigma (0)$ \ and \ $\sigma (1)$. Then the intersection \ $\sigma (I) \cap C$ \ consists of an equal number of positive and negative crossings. As we go along the curve \ $\sigma$ \ there may be consecutive crossings \ $\sigma (t), \sigma (t') \not\in \{ \sigma (0), \sigma (1)\}$ \ with opposite signs. If they lie in the same pathcomponent \ $C'$ \ of \ $C$ \ pick a path in \ $C'$ \ joining them and replace \ $\sigma | (t - \varepsilon, t' + \varepsilon)$ \ by the corresponding \lq\lq parallel\rq\rq\ path just outside of \ $C'$. In the end there must remain consecutive crossings \ $\sigma (t), \sigma (t')$ \ (with opposite signs) which lie in different pathcomponents \ $C'_0$ \ and \ $C'_1$ \ of \ $C$. Apply connected sum surgery along the arc \ $\sigma ([t, t'])$ \ (which can be made embedded via suitable shortcuts). This corresponds to a homotopy of \ $f$ \ which decreases \ $\# \pi_0 (f^{-1} \{*\})$. Iterating this procedure we see that \ $MCC (f, *) \le N (f, *)$ \ if \ $N (f, *) \ne 0$. In view of 1.9 (iii) this completes the proof of theorem 1.13.  \hfill $\blacksquare$

\proclaim{Proposition 6.1} \
Given maps \ $f_1, f_2 : M \to S^1$ \ and a coincidence point \ $x_0 \in M$, put \ $y_0 := f_1 (x_0) = f_2 (x_0)$. Let the kernel of
$$
f_{1 *} \ - \ f_{2 *} \ : \ \pi_1 (M, x_0) \ \longrightarrow \ \pi_1 (S^1, y_0)
$$
act in the standard fashion on the universal covering space \ $(\widetilde M, \widetilde x_0)$ \ of \ $(M, x_0)$ \ by covering transformations.

Then the fiber space \ $E (f_1, f_2)$ \ is fiber homotopically equivalent to
$$
(\widetilde M / \ker (f_{1 *} - f_{2 *})) \ \times \ \coker (f_{1 *} - f_{2 *})
$$
$($where the cokernel has the discrete topology$)$.

If the maps \ $f_1, f_2$ \ are homotopic (e.g.\ when they are coincidence free), then \ $E (f_1, f_2)$ \ is fiber homotopically equivalent to \ $M \times \Bbb Z$.
\endproclaim

\demo{Proof} \ Since \ $S^1$ is a Lie group, \ $E (f_1, f_2)$ \ is fibre homeomorphic to \ $E (f := f_1 \cdot f_2^{-1}, 1)$ \ and hence to the pullback, under $f$, of the starting point fibration of
$$
P (S^1, 1) \ = \ \{ \theta : ([0, 1], 1) \longrightarrow (S^1, 1) \ \text{continuous} \} .
$$
Since every path \ $\theta$ \ can be deformed \ $\rel \{0,1\}$ \ to a path with constant speed, \ $P (S^1, 1)$ \ is fiber homotopically equivalent to the universal covering space of \ $S^1$. Clearly, its pullback under \ $f$ \ allows the indicated description in terms of \ $\widetilde M$, provided a coincidence point exists. This is always the case after a suitable homotopy (cf.\ also 3.2). \hfill
\enddemo

\vskip7mm
\specialhead 7.\ \ The setting of homotopy groups
\endspecialhead

In this section we study the situation where \ $M = S^m$. Moreover, we assume \ $n \ge 2$ \ (having settled the case \ $n = 1$ \ in \S\ 6).

Given basepoint preserving maps \ $f_1, f_2 : (S^m, x_0) \to (N, y_0)$, \ a fixed choice of local orientations of \ $S^m$ \ and \ $N$ \ at \ $x_0$ and $y_0$, resp., determines an isomorphism
$$
i \ : \ \Omega_{m - n}^{fr} (\Lambda (N, y_0)) \ \longrightarrow \ \Omega_{m -n} (E (f_1, f_2); \widetilde\varphi)
\tag 7.1
$$
(apply 2.4 using constant paths). For any bordism class \ $z = [C, \widetilde g, \overline g]$ \ in the target the corresponding class \ $i^{- 1} (z)$ \ can be described as follows. Pick any homotopy \ $G : C \times [0, 1] \to S^m$ \ joining \ $g := pr \scirc \widetilde g$ \ to the constant map at \ $x_0$ \ and lift it in a natural way to a homotopy \ $\widetilde G$ \ in \ $E (f_1, f_2)$ joining \ $\widetilde g$ \ to \ $(x_0, \widetilde g')$, where for \ $x \in C$ the loop \ $\widetilde g' (x)$ \ is the composite of paths
$$
y_0 @>{(f_1 \scirc G (x, -))^{-1}}>> \ f_1 (g (x)) \ @>{\theta (\widetilde g (x))}>> \ f_2 (g (x)) \ @>{f_2 \scirc G (x, -)}>> \ y_0 \ .
\tag 7.2
$$
$\widetilde G$ \ also induces a vector bundle isomorphism joining \ $\overline g$ to a framing \ $\overline g'$ \ of \ $C$. Then
$$
i^{- 1} ([C, \widetilde g, \overline g]) \ = \ [C, \widetilde g', \overline g'] .
\tag 7.3
$$

We use isomorphisms such as \ $i$ \ to identify all \ $\widetilde\omega$-invariants with elements in the {\it same}  \ group which no longer varies with \ $f_1$ or $f_2$. This identification is compatible with base point preserving homotopies (compare 3.2 and 3.3). Thus we obtain well-defined maps 
%% Einschub fuer Grafikeinbindung
$$
% eventuell Breite einstellen
\EPSFxsize=9cm % (andere Angabe)
\BoxedEPSF{diag74.eps}
\tag 7.4
$$
%% Einschub Ende

\proclaim{Proposition 7.5} \
(i) \ For all \ $[f_1], [f_2] \in \pi_m (N, y_0)$
$$\widetilde\omega (f_1, f_2) \ = \ \widetilde\deg (f_1) + (- 1)^n \inv_* (\widetilde\deg (f_2))
$$
where \ $\inv_*$ \ is induced by the involution on \ $\Lambda (N, y_0)$ \ which inverses the direction of a loop;

(ii) \ $\widetilde\deg$ \ is a group homomorphism. 
\endproclaim

\demo{Proof} \ After a base point preserving homotopy and a further small approximation we may assume that \ $f_1$ and $f_2$ \ take different constant values (which lie near \ $y_0$) on complementary halfspheres \ $S^m_{\pm}$, i.e.\
$$
f_1 | S^m_+ \ \equiv \ * \ \ne \ *' \ \equiv \ f_2 | S^m_- .
$$
Then clearly
$$
\widetilde\omega (f_1, f_2) \ = \ \widetilde\omega (f_1, *') \ + \ \widetilde\omega (*, f_2)
$$
and our first claim follows from 4.5 and 7.2.

A very similar argument implies the additivity of the refined degree \ $\widetilde\deg$. \hfill $\blacksquare$
\enddemo

In view of formula 7.5 (i) \ it is natural to expect some interrelation between the following two possible conditions concerning \ $S^m$ and $N$.

\subhead{Condition A}
\endsubhead \ For any two maps \ $f_1, f_2 : S^m \to N$ \ the Nielsen number (or, equivalently, the \ $\widetilde\omega$-invariant) is the only looseness obstruction; in other words: \ $(f_1, f_2)$ \ is loose if and only if \ $N (f_1, f_2) = 0$.

\subhead{Condition B}
\endsubhead \ The refined degree homomorphism \ $\widetilde\deg$ \ in diagram 7.4 is injective.

\proclaim{Proposition 7.6} \ {\rm(i)} Assume that \ $N$ allows a self-homeomorphism a without fixed points $($this holds e.g.\ if \ $N$ \ has a nowhere vanishing vector field$)$. Then condition B implies condition A $($and more specifically that  \ $(f_1, f_2)$ \ is loose if and only if \ $f_2$ \ is homotopic to \ $a \scirc f_1)$.

{\rm(ii)} \ Assume that \ $\pi_m (N-\{\text{point}\}) = 0$. Then condition A implies condition B.
\endproclaim

\demo{Proof} \ (i) \ The assumption implies that for every map \ $f_1 : S^m \to N$ there is a map \ $\widehat f_2$ (e.g.\ $\widehat f_2 = a \circ f_1$) which has no coincidences with \ $f_1$; this is essentially all we need in the proof.

Now consider any pair \ $(f_1, f_2)$ \ with vanishing Nielsen number. After suitable homotopies \ $f_1, f_2$ \ and \ $\widehat f_2$ \ preserve base points. Now \ $\widetilde\omega (f_1, f_2) = \widetilde\omega (f_1, \widehat f_2) = 0$ \ and hence \ $\widetilde\deg (f_2) = \widetilde\deg (\widehat f_2)$ \ by 7.5 (i). If condition B holds \ $f_2$ \ and \ $\widehat f_2$ \ are homotopic and \ $(f_1, f_2)$ \ is loose.

(ii) \ Suppose that condition A holds. Given \ $[f] \in \ker (\widetilde\deg)$, the pair \ $(f_1, f_2) := (f, y_0)$ \ has Nielsen number zero and must be loose. Then, by \cite{Br}, \ $f$ \ is homotopic to a map which avoids \ $y_0$. In view of the assumption in (ii), \ $[f] = 0$ \ and \ $\widetilde\deg$ \ is injective.   \hfill $\blacksquare$

\enddemo 
\vskip7mm
\specialhead 8.\ \ Maps between spheres
\endspecialhead

In this section we work out the details of Example IV of the introduction. Thus let \ $N$ \ be the \ $n$-sphere (with the antipodal selfhomeomorphism denoted by \ $a$). This setting provides a good testing ground for the strength of our methods. Indeed, given {\it any} \ space \ $M$, the answer to the looseness question for maps \ $f_1, f_2 : M \to S^n$ \ is known from the very beginning (cf.\ \cite{DG}, 1.10): if \ $(f_1, f_2)$ \ is homotopic to a coincidence-free pair \ $(f'_1, f'_2)$ \  then \ $f'_1 \sim a \scirc f'_2$ \ via the homotopy \ $((1-t) f'_1 - t f'_2) / \Vert ((1 - t) f'_1 - t f'_2) \Vert$; thus \ $(f_1, f_2)$ \ is loose if and only if \ $f_1 \sim a \scirc f_2$.

Now assume also that \ $M = S^m$ \ and \ $n \ge 2$. 

First we want to make \ $C (f_1, f_2)$ \ pathconnected. After suitable homotopies our maps have the form described in the proof of proposition 7.5, with \ $*'$ \ and \ $*$ \ being regular values of \ $f_1$ \ and \ $f_2$, resp.; moreover \ $f_1$ \ is defined, on a tubular neighbourhood \ $f_1^{-1} \{*'\} \times B^n$ \ of \ $f_1^{-1} \{*'\}$, by the projection to \ $B^n / \partial B^n \cong S^n$ \ (as in the Pontryagin-Thom construction), and similarly \ $f_2$. Then any two pathcomponents of \ $C (f_1, f_2) = f_1^{-1} \{*'\} \cup f_2^{-1} \{*\}$ \ can be joined by an arc \ $\sigma (I)$ \ in \ $S^m$ \ which intersects their tubular neighbourhoods in appropriate rays and such that \ $f_1 (\sigma (I))$ \ and \ $f_2 (\sigma (I))$ \ lie in the same arc from \ $*$ \ to \ $*'$ \ in \ $S^n$. Thus the techniques of the proof of proposition 4.7 allow us to decrease the number of pathcomponents of \ $C (f_1, f_2)$ \ successively (even when \ $n = 2$). This proves the last claim in theorem 1.14.

Next apply the discussion of \ \S\ 7 to the case  \ $N = S^n$.  Then by proposition 7.6 condition A is equivalent to the injectivity of \ $\widetilde\deg$. Most of theorem 1.14 as well as corollary 1.15 follow as soon as we have established the indicated connection with the  James-Hopf invariants.

For this purpose we will exploit the close relationship between loop spaces and immersions. Recall the bijection 
$$
\Cal J_1 (V, \varepsilon^{n -1}) \ @>{\beta \circ \iota}>> \ [S V_c, S^n] \ \cong \ [V_c, \Lambda (S^n, y_0)]
\tag 8.1
$$
discussed in \cite{KS}, theorem 1.2 and example (i) (on p.\ 284). Here \ $V$ \ and \ $(S) V_c$, resp., denote an arbitrary smooth manifold without boundary and (the suspension of) its one-point compactification, resp.. \ $\Cal J_1 (V, \varepsilon^{n -1})$ \ is the set of bordism classes of embeddings \ $e : Q \hookrightarrow V \times \Bbb R$ \ which project to framed \ $(n - 1)$-codimensional immersions in \ $V$; the commuting diagram
%% Einschub fuer Grafikeinbindung
$$
% eventuell Breite einstellen
\EPSFxsize=4cm % (andere Angabe)
\BoxedEPSF{diag82.eps}
\tag 8.2
$$
%% Einschub Ende
subsumes this projection property. By forgetting it we obtain a \ {\it bijection} \ $\iota$ \ onto the bordism set of framed codimension-$n$ embeddings in \ $V \times \Bbb R$. The Pontryagin-Thom procedure then yields a further bijection \ $\beta$ \ onto the indicated base point preserving homotopy set.

Now, given any bordism class in \ $\Omega_*^{fr} (\Lambda (S^n, y_0))$, \ represent it by a singular manifold \ $(V, v)$ \ and apply the inverse bijection \ $(\beta \circ \iota)^{- 1}$ \ (cf.\ 8.1) to the homotopy class of the map \ $v$, to obtain an immersion \ $e' : Q \looparrowright V$ \ as in 8.2.

\proclaim{Proposition 8.3} \
This procedure yields a well-defined isomorphism from the framed bordism group \ $\Omega_r^{fr} (\Lambda (S^n, y_0)), \ r \in \Bbb Z$, onto the joint bordism group \ $I_{r, n -1}$ \ of quadruples \ $(e', Q, V, e'')$ \ where
\roster
\item"(i)" \ $e' : Q \looparrowright V$ \ is a framed \ $(n -1)$-codimensional immersion between closed smooth manifolds;
\item"(ii)" \ $V$ \ is framed and has dimension \ $r$; \ and
\item"(iii)" \ $e'' : Q \to \Bbb R$, together with \ $e'$, determines an embedding \ $e = (e', e'')$ \ (which \lq\lq decompresses\rq\rq \ $e'$) as in $8.2$.
\endroster
\endproclaim

A \ {\it joint bordism} \ of such quadruples is free to alter \ $Q$ and $V$ \ simultaneously (compare \cite{S}).

\demo{Proof} \ Given a bordism \ $(W, w)$ \ between, say, \ $(V_1, v_1)$ \ and \ $(V_2, v_2)$, we can replace \ $(V, v)$ \ in the above-mentioned construction by the closed singular manifold \ $(W \cup_\partial - W, w \cup w)$ \ and obtain the required joint bordism between the decompressed immersions which correspond to  \ $(V_1, v_1)$ \ and \ $(V_2, v_2)$. Thus our construction yields a well-defined inverse of the obvious homomorphism induced by such forgetful maps as \ $\iota$ \ in 8.1. \hfill $\blacksquare$
\enddemo

Next, given a quadruple as in 8.3, we analyze the multiple self-intersections of \ $e'$. After a small deformation we may assume that they are all transverse. Then for \ $k \ge 1$ \ the set of \ $k$-tuple points
$$
Q_k  :=  \{ (x_1, . . , x_k) \in Q^k | e' (x_1) = . . = e' (x_k), e'' (x_1) < e'' (x_2) < . . < e'' (x_k) \}
\tag 8.4
$$
is a closed manifold, smoothly immersed in \ $V$ \ with a natural framing (inherited from the framing of \ $e'$ \ and the ordering of the intersection branches given by \ $e''$). Using the stable parallelization of \ $V$ \ and forgetting the immersion, we obtain the stable bordism class
$$
h_{k + 1} ([e', Q, V, e'']) \ := \ [Q_k] \ \ \in \ \ \Omega^{fr}_{r -k (n -1)} \cong \pi^S_{r -k (n - 1)}
\tag 8.5
$$
(note the shift of the index; e.g.\ the \ $h_2$-value equals \ $[Q]$). Also we put  
$$
h_1 ([e', Q, V, e'']) := [V] \ \ .
$$

\proclaim{Proposition 8.6} \
This construction determines a well-defined isomorphism
$$
h = \bigoplus_{k \ge 1} \ h_k \ : \ I_{r, n -1} \ @>{ \ \cong \ }>> \ \bigoplus_{k \ge 1}  \pi^S_{r - (k -1)(n -1)}
$$
$($compare $8.3)$.
\endproclaim

\demo{Proof} \ Joint bordisms (and disjoint unions) yield bordisms (and unions) of self-intersections; thus \ $h$ \ is a well-defined homomorphism.

As for bijectivity, we follow an unpublished argument of Paul Schweitzer (compare also \cite{K 2}, pp.\ 81--83). Given \ $[e', Q, V, e''] \in \ker h$, we may assume that \ $e'$ \ is {\it self-transverse}, i.e.\ all self-intersections are transverse. Let \ $k$ \ be the highest occuring multiplicity. Then \ $Q_k$ \ is an embedded submanifold of \ $V$ \ with a tubular neighbourhood of the form \ $Q_k \times (\Bbb R^{n -1})^k$ \ which intersects \ $e' (Q)$ \ in the branches \ $Q_k \times (\Bbb R^{n -1})^j \times \{0\} \times (\Bbb R^{n -1})^{k -j -1}, \ j = 0, \dots, k -1$.  Now pick a nulbordism of \ $Q_k$ \ and glue its product with the unit ball \ $B_k$ \ of \ $(\Bbb R^{n -1})^k$ \ to \ $V \times [0, 1]$ \ along \ $Q_k \times B_k \subset V \times \{1\}$. The resulting joint bordism \lq\lq seals off\rq\rq \ $Q_k$ \ so that only self-intersections of strictly smaller multiplicity survive. Iterating this procedure we obtain a nulbordism of \ $(e', Q, V, e'')$. Thus \ $h$ \ is injective.

Surjectivity follows in a similar way. Given a closed framed manifold \ $Q'_k$ \ of dimension \ $r - k (n -1)$, equip the product \ $Q'_k \times B_k$ \ with the obvious codimension-$(n -1)$ immersion as described above, seal off the self-intersections of multiplicity \ $k -1$ \ in the boundary etc.\ to produce \ $z \in I_{r, n -1}$ \ such that \ $h (z) = [Q'_k] + h_\ell$--values for \ \ $\ell < k +1$.

It is not hard to extend \ $e''$ in these constructions since they preserve the ordering of the branches of the immersion \ $e'$ \ at the self-intersections. \ \hfill $\blacksquare$
\enddemo

Finally we turn to the diagram in theorem 1.14. According to propositions 8.3 and 8.6 we may identify \ $\Omega^{fr}_{m -n} (\Lambda (S^n, y_0))$ \ with \ $I_{m -n, n -1}$ \ and then apply the selfintersection isomorphism \ $h$.

Let us calculate \ $h \scirc \widetilde\deg$. Given \ $[f] \in \pi_m (S^n, y_0)$, pick a regular value \ $\ast \in S^n$ \ near, but different from, the basepoint \ $y_0 = f (x_0)$. In view of the bijection 8.1, applied to \ $V = \Bbb R^{m -1}$, we may assume that the inclusion
$$
e_C \ : \ C = f^{- 1} (\{\ast\}) \ \ \ \subset \ S^m - \{ x_0 \} \ = \ \Bbb R^m
\tag 8.7
$$
of the coincidence manifold projects to a framed immersion \ $e'_C$ \  into \ $\Bbb R^{m -1}$. Now let \ $G$ \ denote a contraction of \ $e_C$ \ in the negative \ $x_m$--direction towards \ $\infty = x_0 \in S^m$. Then $f \circ G$, together with a path \ $\theta_*$ \ in \ $S^n$ \ from \ $\ast$ \ to \ $y_0$, describe a map
$$
v \ : \ (C  \times [0, 1], C \times \{0, 1\} ) \ \longrightarrow \ (S^n, y_0)
$$
such that
$$
\widetilde\deg (f) \ = \ [C, \ \text{adjoint of}\ \ v]
\tag 8.8
$$
(compare 3.2, 7.1, and 7.2).

In order to determine \ $h \circ \widetilde\deg (f)$, \ we now apply the discussion of 8.1 to the case \ $V = C$ \ and to the map \ $v$. \ Pick a point \ $\ast' \in S^n - \{ y_o\} \ = \ \Bbb R^n$ \ near \ $\ast$, \ but with strictly larger \ $x_n$--coordinate and such that \ $\ast'$ \ also avoids the image of the path \ $\theta_*$. \ Then the inclusion
$$
e_Q \ = \ (e'_Q, e''_Q) \ \ : \ \ Q := v^{- 1} (\{ *'\} ) \ \ \ \subset \ \ \ C \times \Bbb R
$$
projects to a framed immersion \ $e'_Q$ \ of codimension \ $n -1$; the branches of \ $e'_Q$ \ through a point \ $c \in C$ \ correspond to those branches of the immersion \ $e'_C \ : \ C \ \looparrowright \ \Bbb R^{m -1}$ \ (compare 8.7) which pass \ {\it under} \ the point \ $c$ \ w.r.\ to the \ $x_m$--coordinate \ $e''_C$. \ Therefore we can identify \ $\widetilde\deg (f)$ \ (cf.\ 8.8), via the isomorphism in proposition 8.3, with the joint bordism class \ $[e'_Q, Q, C, e''_Q] \in I_{m-n, n-1}$. Moreover the \ $(k -1)$--tuple point manifold of \ $e'_Q$ \ (which represents \ $h_k (\widetilde\deg (f))$, cf.\ 8.5) equals essentially the $k$-tuple point locus of the \lq\lq compression\rq\rq \ $e'_C$ \ of \ $e_C$ \ (cf. 8.7). Thus, according to \cite{KS} (see, in particular, p.\ 287, l.\ 20, and theorem 3.2) it also represents the stabilized James-Hopf invariant (cf.\ \cite{J}) -- at least up to a sign which depends only on \ $k$ and which we subsume into the definition of \ $h_k$.

In view of proposition 7.5 and of the identity
$$
0 \ = \ \widetilde\omega (a \scirc f, f) \ = \ \widetilde\deg (a \scirc f) + (- 1)^n \inv_* (\widetilde\deg (f))
$$
the remaining claims in theorem 1.14 follow from

\proclaim{Proposition 8.9} \ Let \ $\inv$ \ denote the involution on \ $\Lambda (S^n, y_0)$ \ which inverses the direction of a loop, and let \ $\inv_*$ \ be the induced involution on \ $\Omega_*^{fr} (\Lambda (S^n, y_0))$.

Then for \ $k \ge 1$
$$
\qquad \qquad h_k \ \scirc \ \inv_* \ \ \ \ \ \ = \ \ \ \ \ \ - (- 1)^{k + (n -1) {k - 1 \choose 2}} h_k  \ \ ;
$$
moreover for all \ $[f] \in \pi_* (S^n)$ \ we have
$$
\Gamma_k (f) := E^\infty \scirc \gamma_k (f) \  \ = \ \ \ - (- 1)^{k + (n -1){k \choose 2}} \Gamma_k (f)
$$
and
$$
\Gamma_k (r \scirc f) \ \ \ \ \ = \ \ \ \ \ (- 1)^k \Gamma_k (f)
$$
where \ $r : (S^n, y_0) \to (S^n; y_0)$ \ is a reflection. (Note that \ $a \sim r^{n +1})$.

\endproclaim

\demo{Proof} \ If in the construction of \ $h_k$ a given element \ $\lambda \in \Omega^{fr}_* (\Lambda (S^n, y_0))$ \ corresponds to a quadruple \ $(e', Q, V, e'')$ \ as in 8.2 and 8.3 \ then \ $\inv_* (\lambda)$ \ corresponds to the \lq\lq reflected\rq\rq\ quadruple \ $(\overline e', Q, V, -e'')$. Here \ $e'$ \ and \ $\overline e'$ \ denote the same immersion, but with opposite framing: one of the frame vectors is replaced by its negative; this is needed to make up for the reflection in the \ $\Bbb R$-component which takes \ $e''$ to $ -e''$. \ Thus when we compare \ $h_k (\inv_* (\lambda))$ \ to \ $h_k (\lambda)$ \ each intersection branch at the \ $(k -1)$-tuple point set \ $Q_{k -1}$ (cf.\ 8.4) contributes a factor \ $- 1$. Moreover the reflection reverses the ordering of the \ $(n -1)$-codimensional branches; interchanging them yields \ $(n -1)((k -2) + (k -1) + \dots + 2 + 1)$ sign changes. All this adds up to the factor \ $(- 1)^{k -1 + (n -1){k -1 \choose 2}}$ \ claimed in proposition 8.9.

A similar argument shows that the diagram
$$
\CD
\pi_m (S^n) \ \cong \ \  @. \pi_{m -1} (\Lambda S^n)   @>{ \ \ \Gamma_k \ \ }>>  \pi^S_{m -1 - k (n -1)} \\
@VV{- \id}V                     @VV{\inv_{*}}V                   @VV{(- 1)^{k + (n -1){k \choose 2}} \cdot \id}V \\
\pi_m (S^n) \ \cong \ \  @. \pi_{m -1} (\Lambda S^n)   @>{ \ \ \Gamma_k \ \ }>>  \pi^S_{m -1 - k (n -1)}
\endCD
$$

\flushpar
commutes. Also any given reflection \ $r$ of $S^n$ \ can be deformed until it preserves the splitting \ \ $S^n - \{ \infty\} \ \cong \ \Bbb R^{n -1} \times \Bbb R$ \ \ together with the \ $\Bbb R$-values. Thus it preserves also the ordering of the selfintersection branches, but it reverses the normal framing of a projected immersion \ $e'$ \ into \ $V = \Bbb R^{m -1}$ \ (compare 8.2). At the \ $k$-tuple point set this amounts to a \ $k$-fold change of signs. \hfill $\blacksquare$
\enddemo

\example{Example 8.10} \ The image of the homomorphism
$$
\widetilde\deg_2 \ = \ \Gamma_2 \ \ : \ \ \pi_{12} (S^5) \cong \Bbb Z_{30} \ \longrightarrow \ \pi^S_3 \cong \Bbb Z_{24}
$$
lies in a group of order 2. 
\endexample

\bigskip
\demo{Proof of corollary {\rm 1.16}} \
We use the tables in chapter XIV of Toda's book \cite{T}. Also we identify elements of \ $\pi_m (S^n)$ \ via the Pontryagin-Thom  construction with bordism classes of closed smoothly embedded \ $(m - n)$-manifolds \ $Q \subset \Bbb R^m$ \ which are equipped with a framing (i.e.\ a nonstable trivialisation of the normal bundle \ $\nu (Q, \Bbb R^m))$. The stable suspension homomorphism \ $E^\infty$ \ forgets the embedding and retains only the stably framed manifold \ $Q$.

Thus we can represent e.g.\ the generator \ $\eta_2$ \ of \ $\pi_3 (S^2) \cong \Bbb Z$ \ by a framed circle in \ $\Bbb R^3$ \ which projects to the framed figure-8 immersion in \ $\Bbb R^2$, and we see that
$$
\Gamma (\eta_2) \ = \ (1, 1) \ \in \ \Bbb Z_2 \oplus \Bbb Z \ \cong \ \pi^S_1 \oplus \pi^S_0 \ .
$$
Similarly \ $\eta_2 \scirc \eta_3 \in \pi_4 (S^2)$ \ and \ $\eta_2 \scirc \eta_3 \scirc \eta_4 \in \pi_5 (S^2)$ \ are represented by embedded tori \ $(S^1)^2$ \ and \ $(S^1)^3$ \ which suspend to the nontrivial elements \ $\eta^2 \in \pi^S_2$ \ and \ $\eta^3 \in \pi^S_3$, resp. It follows that \ $\Gamma$ \ is injective on \ $\pi_{n +1} (S^n), \ n \ge 2$, as well as on the groups \ $\pi_{n + 2} (S^n), \ n \ge 2$, and \ $\pi_5 (S^2)$ \ which all have order 2.

To finish the proof consider the following exact pieces of EHP-sequences (cf. \cite{W}, p. 542)
$$
% eventuell Breite einstellen
\EPSFxsize=10cm % (andere Angabe)
\BoxedEPSF{diag2.eps}
$$
and
$$
% eventuell Breite einstellen
\EPSFxsize=9cm % (andere Angabe)
\BoxedEPSF{diag3.eps}
$$
%% Einschub Ende
\medskip

Clearly \ $H = E^\infty \scirc \gamma_2$ \ is onto here; moreover, \ $E$ \ and  \ $E^\infty \oplus H$ \ are injective. This implies our claim also for \ $\pi_6 (S^3)$ \ and \ $\pi_7 (S^4)$. \hfill $\blacksquare$
\enddemo

\demo{Proof of corollary {\rm 1.17}} \
First recall that for odd \ $n$ \ the Whitehead product \ $[f] := [\iota_n, \iota_n]$ \ lies in the kernel of the suspension homomorphism \ $E$ (and hence of \ $E^\infty$) and also of \ $\gamma_2$ \ (since \ $2 [\iota_n, \iota_n] = 0$, cf.\ \cite{W}, p.\ 474 and 485). Moreover, by the famous result of F.\ Adams on odd Hopf invariants (and again by an EHP-sequence argument) \ $[\iota_n, \iota_n] \ne 0$ \ if \ $n \ne 1, 3, 7$.

In the remaining cases of corollary 1.17 an inspection of Toda's tables (cf. \cite{T}, p.\ 186--188) shows that \ $\Gamma$ \ is also not injective: the domain turns out to be larger than the relevant torsion part of the target group, or the orders of group elements are not compatible; in view of the relation \ $2 \Gamma_4 \equiv 0$ \ this argument works also for
$$
\Gamma \ : \ \pi_{24} (S^6) \ = \ \Bbb Z_{24} \oplus \Bbb Z_6 \oplus \Bbb Z_2 \longrightarrow (\Bbb Z_8 \oplus \Bbb Z_2) \oplus (\Bbb Z_6 \oplus \Bbb Z_2) \oplus (\Bbb Z_2 \oplus \Bbb Z_2) \oplus \Bbb Z_{24} .
$$ \hfill $\blacksquare$

\enddemo

\vskip1cm

\head{References}
\endhead

\enddocument